\documentclass[11pt]{article}
\usepackage{amssymb}
\usepackage{stmaryrd}
\usepackage{indentfirst}
\usepackage{amscd}
\usepackage{amsmath}
\usepackage{latexsym}
\usepackage{amsfonts}
\usepackage{amssymb}
\usepackage{amsthm}
\usepackage{eqnarray}
\usepackage{mathrsfs}
\usepackage{ams,amsmath,amsthm,amsfonts,amssymb,bm}
\usepackage[dvips]{graphicx}
\usepackage{subfigure}
\usepackage{color}

\begin{document}
\title{A sufficient condition on operator order\\ for strictly
positive operators
 \thanks
{ This work is supported by National Natural Science Fund of China
(11171013 ).}}
\author{ Jian Shi$^{1}$ \ { and} \ Zongsheng Gao \hskip 1cm  \\
{\small   LMIB $\&$ School of Mathematics and Systems Science,}\\
{\small Beihang University, Beijing, 100191, China}\\ }
         \date{}
         \maketitle

\maketitle \baselineskip 16pt \line(1,0){420}

\noindent{\bf Abstract.}\ \  Let $A_{1}$, $A_{2}$, $\ldots$, $A_{k}$
be strictly positive operators on a Hilbert space. This note is to
show a sufficient condition of $A_{k}\geq A_{k-1}\geq \cdots\geq
A_{3}\geq A_{2}\geq A_{1}$, which extends the related result
before.\vspace{0.2cm}


\vspace{0.2cm} \vskip   0.2cm \footnotetext[1]{Corresponding author.

E-mail addresses: shijian@ss.buaa.edu.cn\ (J. Shi),
zshgao@buaa.edu.cn\ (Z. Gao).}

\setlength{\baselineskip}{20pt}
\section{Introduction }

In this note, we denote a bounded linear operator on a Hilbert space
by a capital letter, such as $T$. $T\geq 0$ and $T>0$ stand for a
positive operator and a strictly positive operator, respectively. In
what follows, we assume that $A_{1}$, $A_{2}$, $\ldots$, $A_{k}$ are
strictly positive operators.

In \cite{Lin_2007}, C. -S. Lin proved several characterizations of
operator order $A_{2}\geq A_{1}$ in terms of Furuta
inequality\cite{Furuta_1987} and Pedersen-Takesaki type operator
equation\cite{Pedersen-Takesaki}; Afterwards, C. -S. Lin and Y. J.
Cho showed characterizations of $A_{3}\geq A_{2}\geq A_{1}$ in
\cite{Lin_2011} by extended grand Furuta
inequality\cite{Uchiyama_2003}; As generalizations, J. Shi and Z.
Gao gave characterizations of $A_{k}\geq A_{k-1}\geq \cdots\geq
A_{3}\geq A_{2}\geq A_{1}$ in\cite{Shi} by Further extension of the
grand Furuta inequality\cite{Yang_2010}. As a continuation, this
note is to prove a sufficient condition of  $A_{k}\geq A_{k-1}\geq
\cdots\geq A_{3}\geq A_{2}\geq A_{1}$.

Let us recall two important theorems first.

\noindent {\bf Theorem 1.1} (L\"{o}wner-Heinz inequality)
(\cite{Heinz, Lowner}). If $P\geq Q\geq 0$, then $P^{\alpha}\geq
Q^{\alpha}$ holds for any $\alpha\in$ [0, 1].

\noindent {\bf Theorem 1.2} (\cite{Yuan_Wang}). For $P, Q>0$,
$r+\delta>0$ with $r>0$, $w\in [0, 1]$. If $P^{r+\delta}\geq
(P^{\frac r 2}Q^{s}P^{\frac r 2})^{w}$ holds for any $s>1$, then
$Q\leq I$; If $P^{r+\delta}\leq (P^{\frac r 2}Q^{s}P^{\frac r
2})^{w}$ holds for any $s>1$, then $Q\geq I$.

\section{A sufficient condition of $A_{k}\geq A_{k-1}\geq \cdots\geq
A_{3}\geq A_{2}\geq A_{1}$}

This section is to show a sufficient condition of $A_{k}\geq
A_{k-1}\geq \cdots\geq A_{3}\geq A_{2}\geq A_{1}$. First, we
consider the condition that $k$ is an odd integer.

\noindent {\bf Theorem 2.1.} For $t_{1}, t_{2}, \ldots, t_{n},
w_{1}, w_{2}, \ldots, w_{2n}\in [0, 1]$, $r>t_{n}$. If the following
inequalities always hold for $p_{1}, p_{2}, \ldots, p_{2n}\geq 1$,\\
\indent (I.1) $A^{r-t_{n}}_{2n+1}\geq     \Big\{A^{\frac r 2}_{2n+1}\Big[A^{-{\frac {t_{n}}{2}}}_{2n}\big\{A^{\frac {t_{n-1}}{2}}_{2n-1} \cdots A^{\frac {t_{2}}{2}}_{5}\big[A^{-{\frac {t_{2}}{2}}}_{4}\cdot\{A^{\frac {t_{1}}{2}}_{3}(A^{-{\frac {t_{1}}{2}}}_{2}A^{p_{1}}_{1}A^{-{\frac {t_{1}}{2}}}_{2})^{p_{2}} A^{\frac {t_{1}}{2}}_{3}\}^{p_{3}}\cdot \\
 A^{-{\frac {t_{2}}{2}}}_{4}\big]^{p_{4}}A^{\frac {t_{2}}{2}}_{5}
\cdots A^{\frac {t_{n-1}}{2}}_{2n-1}\big\}^{p_{2n-1}}A^{-{\frac {t_{n}}{2}}}_{2n}\Big]^{p_{2n}}A^{\frac r 2}_{2n+1}\Big\}^{w_{1}}$;\\
\indent (I.2) $A^{r-t_{n}}_{2n+1}\geq     \Big\{A^{\frac r 2}_{2n+1}\Big[A^{-{\frac {t_{n}}{2}}}_{2n+1}\big\{A^{\frac {t_{n-1}}{2}}_{2n } \cdots A^{\frac {t_{2}}{2}}_{6}\big[A^{-{\frac {t_{2}}{2}}}_{5}\cdot\{A^{\frac {t_{1}}{2}}_{4}(A^{-{\frac {t_{1}}{2}}}_{3}A^{p_{1}}_{2}A^{-{\frac {t_{1}}{2}}}_{3})^{p_{2}} A^{\frac {t_{1}}{2}}_{4}\}^{p_{3}}\cdot \\
 A^{-{\frac {t_{2}}{2}}}_{5}\big]^{p_{4}}A^{\frac {t_{2}}{2}}_{6}
\cdots A^{\frac {t_{n-1}}{2}}_{2n }\big\}^{p_{2n-1}}A^{-{\frac {t_{n}}{2}}}_{2n+1}\Big]^{p_{2n}}A^{\frac r 2}_{2n+1}\Big\}^{w_{2}}$;\\
\indent (I.3) $A^{r-t_{n}}_{2n+1}\geq     \Big\{A^{\frac r 2}_{2n+1}\Big[A^{-{\frac {t_{n}}{2}}}_{2n+1}\big\{A^{\frac {t_{n-1}}{2}}_{2n+1 } \cdots A^{\frac {t_{2}}{2}}_{7}\big[A^{-{\frac {t_{2}}{2}}}_{6}\cdot\{A^{\frac {t_{1}}{2}}_{5}(A^{-{\frac {t_{1}}{2}}}_{4}A^{p_{1}}_{3}A^{-{\frac {t_{1}}{2}}}_{4})^{p_{2}} A^{\frac {t_{1}}{2}}_{5}\}^{p_{3}} \cdot\\
 A^{-{\frac {t_{2}}{2}}}_{6}\big]^{p_{4}}A^{\frac {t_{2}}{2}}_{7}
\cdots A^{\frac {t_{n-1}}{2}}_{2n+1 }\big\}^{p_{2n-1}}A^{-{\frac {t_{n}}{2}}}_{2n+1}\Big]^{p_{2n}}A^{\frac r 2}_{2n+1}\Big\}^{w_{3}}$;\\
\indent  {  $\cdots \cdots \cdots \cdots$}\\
\indent (I.n) $A^{r-t_{n}}_{2n+1}\geq     \Big\{A^{\frac r 2}_{2n+1}\Big[A^{-{\frac {t_{n}}{2}}}_{2n+1}\big\{A^{\frac {t_{n-1}}{2}}_{2n+1 } \cdots A^{\frac {t_{2}}{2}}_{n+4}\big[A^{-{\frac {t_{2}}{2}}}_{n+3}\{A^{\frac {t_{1}}{2}}_{n+2}(A^{-{\frac {t_{1}}{2}}}_{n+1}A^{p_{1}}_{n}A^{-{\frac {t_{1}}{2}}}_{n+1})^{p_{2}} A^{\frac {t_{1}}{2}}_{n+2}\}^{p_{3}}  \\
 A^{-{\frac {t_{2}}{2}}}_{n+3}\big]^{p_{4}}A^{\frac {t_{2}}{2}}_{n+4}
\cdots A^{\frac {t_{n-1}}{2}}_{2n+1 }\big\}^{p_{2n-1}}A^{-{\frac {t_{n}}{2}}}_{2n+1}\Big]^{p_{2n}}A^{\frac r 2}_{2n+1}\Big\}^{w_{n}}$;\\
\indent (I.n+1) $A^{r-t_{n}}_{ 1}\leq     \Big\{A^{\frac r 2}_{ 1}\Big[A^{-{\frac {t_{n}}{2}}}_{1}\big\{A^{\frac {t_{n-1}}{2}}_{ 1 } \cdots A^{\frac {t_{2}}{2}}_{n-2}\big[A^{-{\frac {t_{2}}{2}}}_{n-1}\{A^{\frac {t_{1}}{2}}_{n }(A^{-{\frac {t_{1}}{2}}}_{n+1}A^{p_{1}}_{n+2}A^{-{\frac {t_{1}}{2}}}_{n+1})^{p_{2}} A^{\frac {t_{1}}{2}}_{n }\}^{p_{3}}  \\
 A^{-{\frac {t_{2}}{2}}}_{n-1}\big]^{p_{4}}A^{\frac {t_{2}}{2}}_{n-2}
\cdots A^{\frac {t_{n-1}}{2}}_{1}\big\}^{p_{2n-1}}A^{-{\frac {t_{n}}{2}}}_{1}\Big]^{p_{2n}}A^{\frac r 2}_{1}\Big\}^{w_{n+1}}$;\\
\indent  {  $\cdots \cdots \cdots \cdots$}\\
\indent (I.2n-2) $A^{r-t_{n}}_{ 1}\leq     \Big\{A^{\frac r 2}_{
1}\Big[A^{-{\frac {t_{n}}{2}}}_{1}\big\{A^{\frac {t_{n-1}}{2}}_{1}
\cdots A^{\frac {t_{2}}{2}}_{2n-5}\big[A^{-{\frac
{t_{2}}{2}}}_{2n-4}\{A^{\frac {t_{1}}{2}}_{2n-3 }\cdot(A^{-{\frac
{t_{1}}{2}}}_{2n-2}A^{p_{1}}_{2n-1}A^{-{\frac
{t_{1}}{2}}}_{2n-2})^{p_{2}}\cdot \\ A^{\frac
{t_{1}}{2}}_{2n-3}\}^{p_{3}}
 A^{-{\frac {t_{2}}{2}}}_{2n-4}\big]^{p_{4}}A^{\frac {t_{2}}{2}}_{2n-5}
\cdots A^{\frac {t_{n-1}}{2}}_{1}\big\}^{p_{2n-1}}A^{-{\frac {t_{n}}{2}}}_{1}\Big]^{p_{2n}}A^{\frac r 2}_{1}\Big\}^{w_{2n-2}}$;\\
\indent (I.2n-1) $A^{r-t_{n}}_{ 1}\leq     \Big\{A^{\frac r 2}_{
1}\Big[A^{-{\frac {t_{n}}{2}}}_{1}\big\{A^{\frac {t_{n-1}}{2}}_{2}
\cdots A^{\frac {t_{2}}{2}}_{2n-4}\big[A^{-{\frac
{t_{2}}{2}}}_{2n-3}\{A^{\frac {t_{1}}{2}}_{2n-2 }\cdot(A^{-{\frac
{t_{1}}{2}}}_{2n-1}A^{p_{1}}_{2n}A^{-{\frac
{t_{1}}{2}}}_{2n-1})^{p_{2}}\cdot \\ A^{\frac
{t_{1}}{2}}_{2n-2}\}^{p_{3}}
 A^{-{\frac {t_{2}}{2}}}_{2n-3}\big]^{p_{4}}A^{\frac {t_{2}}{2}}_{2n-4}
\cdots A^{\frac {t_{n-1}}{2}}_{2}\big\}^{p_{2n-1}}A^{-{\frac {t_{n}}{2}}}_{1}\Big]^{p_{2n}}A^{\frac r 2}_{1}\Big\}^{w_{2n-1}}$;\\
\indent (I.2n) $A^{r-t_{n}}_{ 1}\leq     \Big\{A^{\frac r 2}_{
1}\Big[A^{-{\frac {t_{n}}{2}}}_{2}\big\{A^{\frac {t_{n-1}}{2}}_{3}
\cdots A^{\frac {t_{2}}{2}}_{2n-3}\big[A^{-{\frac
{t_{2}}{2}}}_{2n-2}\{A^{\frac {t_{1}}{2}}_{2n-1}\cdot(A^{-{\frac
{t_{1}}{2}}}_{2n}A^{p_{1}}_{2n+1}A^{-{\frac
{t_{1}}{2}}}_{2n})^{p_{2}}\cdot \\ A^{\frac
{t_{1}}{2}}_{2n-1}\}^{p_{3}}
 A^{-{\frac {t_{2}}{2}}}_{2n-2}\big]^{p_{4}}A^{\frac {t_{2}}{2}}_{2n-3}
\cdots A^{\frac {t_{n-1}}{2}}_{3}\big\}^{p_{2n-1}}A^{-{\frac
{t_{n}}{2}}}_{2}\Big]^{p_{2n}}A^{\frac r 2}_{1}\Big\}^{w_{2n}}$,\\
then the operator order $A_{2n+1}\geq A_{2n}\geq
A_{2n-1}\geq\cdots\geq A_{3}\geq A_{2}\geq A_{1}$ holds.

\noindent {\bf Proof.} Applying Theorem 1.2 to (I.1), we have
\begin{equation}\tag{2.1}
\begin{split}
&\ A^{-{\frac {t_{n}}{2}}}_{2n}\big\{A^{\frac {t_{n-1}}{2}}_{2n-1}
\cdots A^{\frac {t_{2}}{2}}_{5}\big[A^{-{\frac
{t_{2}}{2}}}_{4}\{A^{\frac {t_{1}}{2}}_{3}(A^{-{\frac
{t_{1}}{2}}}_{2}A^{p_{1}}_{1}A^{-{\frac {t_{1}}{2}}}_{2})^{p_{2}}
A^{\frac {t_{1}}{2}}_{3}\}^{p_{3}}
 A^{-{\frac {t_{2}}{2}}}_{4}\big]^{p_{4}}A^{\frac {t_{2}}{2}}_{5}
\cdots A^{\frac {t_{n-1}}{2}}_{2n-1}\big\}^{p_{2n-1}}A^{-{\frac
{t_{n}}{2}}}_{2n}\\
\leq &I.
\end{split} \end{equation}

\noindent Using L\"{o}wner-Heinz inequality to (2.1) for 2n-2 times,
the following result hold.
\begin{equation}\tag{2.2}
\begin{split}
& \ A_{2}^{-{\frac {t_{1}}{2}}} A_{1}^{p_{1}}A_{2}^{-{\frac
{t_{1}}{2}}}\\
\leq & \Big\{A_{3}^{-{\frac {t_{1}}{2}}}\big[A_{4}^{{\frac
{t_{2}}{2}}}\big(\cdots(A_{2n-1}^{-{\frac
{t_{n-1}}{2}}}A_{2n}^{{\frac {t_{n}}{p_{2n-1}}}}A_{2n-1}^{-{\frac
{t_{n-1}}{2}}})^{\frac {1}{p_{2n-2}}}\cdots\big)^{\frac
{1}{p_{4}}}A_{4}^{{\frac {t_{2}}{2}}}\big]^{\frac
{1}{p_{3}}}A_{3}^{-{\frac {t_{1}}{2}}}\Big\}^{\frac {1}{p_{2}}}.
\end{split}
\end{equation}
Because each $A_{i}$ is strictly positive, there exist a positive
constant $\delta_{i}$ such that $A_{i}\geq
{\frac{1}{\delta_{i}}}I>0$. Therefore,
\begin{equation}\tag{2.3}
\begin{split}
A_{2}^{-{\frac {t_{1}}{2}}} A_{1}^{p_{1}}A_{2}^{-{\frac
{t_{1}}{2}}}\leq \big\{[(\|A_{2n}\|^{\frac
{t_{n}}{p_{2n-1}}}\delta_{2n-1}^{ {t_{n-1}}})^{\frac
{1}{p_{2n-2}}}\cdots\|A_{4}\|^{{t_{2}}}]^{\frac
{1}{p_{3}}}\delta_{3}^{{t_{1}}}\big\}^{\frac {1}{p_{2}}}.
\end{split}
\end{equation}
$A_{2}\geq A_{1}$ holds by putting $t_{1}=p_{1}=1$ and
$p_{2}\rightarrow \infty$ above.

Similarly, we can obtain $A_{3}\geq A_{2}$ by (I.2), $A_{4}\geq
A_{3}$ by (I.3), $\cdots $, $A_{2n+1}\geq A_{2n}$ by (I.2n),
respectively.  $\square$

\noindent {\bf Remark 2.1.} If $w_{1}=w_{2}=\cdots=w_{2n}={\frac
{r-t_{n}}{\psi [2n]-t_{n}+r}}$, where $\psi
[2n]=\{\cdots[\{[(p_{1}-t_{1})p_{2}+t_{1}]p_{3}-t_{2}\}p_{4}+t_{2}]p_{5}-\cdots
-t_{n}\}p_{2n}+t_{n}$, the condition in Theorem 2.1 is the
sufficient and necessary condition of $A_{2n+1}\geq A_{2n}\geq
A_{2n-1}\geq \cdots\geq A_{3}\geq A_{2}\geq A_{1}$. See \cite{Shi}
for details.

Next, we consider that the condition that $k$ is an even integer.\\

\noindent {\bf Theorem 2.2.} For $t_{1}, t_{2}, \ldots, t_{n},
w_{1}, w_{2}, \ldots, w_{2n-1}\in [0, 1]$, $r>t_{n}$. If the
following
inequalities always hold for $p_{1}, p_{2}, \ldots, p_{2n}\geq 1$,\\
\indent (II.1) $A^{r-t_{n}}_{2n}\geq     \Big\{A^{\frac r 2}_{2n }\Big[A^{-{\frac {t_{n}}{2}}}_{2n}\big\{A^{\frac {t_{n-1}}{2}}_{2n-1} \cdots A^{\frac {t_{2}}{2}}_{5}\big[A^{-{\frac {t_{2}}{2}}}_{4}\cdot\{A^{\frac {t_{1}}{2}}_{3}(A^{-{\frac {t_{1}}{2}}}_{2}A^{p_{1}}_{1}A^{-{\frac {t_{1}}{2}}}_{2})^{p_{2}} A^{\frac {t_{1}}{2}}_{3}\}^{p_{3}}\cdot \\
 A^{-{\frac {t_{2}}{2}}}_{4}\big]^{p_{4}}A^{\frac {t_{2}}{2}}_{5}
\cdots A^{\frac {t_{n-1}}{2}}_{2n-1}\big\}^{p_{2n-1}}A^{-{\frac {t_{n}}{2}}}_{2n}\Big]^{p_{2n}}A^{\frac r 2}_{2n }\Big\}^{w_{1}}$;\\
\indent (II.2) $A^{r-t_{n}}_{2n }\geq     \Big\{A^{\frac r 2}_{2n }\Big[A^{-{\frac {t_{n}}{2}}}_{2n }\big\{A^{\frac {t_{n-1}}{2}}_{2n } \cdots A^{\frac {t_{2}}{2}}_{6}\big[A^{-{\frac {t_{2}}{2}}}_{5}\cdot\{A^{\frac {t_{1}}{2}}_{4}(A^{-{\frac {t_{1}}{2}}}_{3}A^{p_{1}}_{2}A^{-{\frac {t_{1}}{2}}}_{3})^{p_{2}} A^{\frac {t_{1}}{2}}_{4}\}^{p_{3}}\cdot \\
 A^{-{\frac {t_{2}}{2}}}_{5}\big]^{p_{4}}A^{\frac {t_{2}}{2}}_{6}
\cdots A^{\frac {t_{n-1}}{2}}_{2n }\big\}^{p_{2n-1}}A^{-{\frac {t_{n}}{2}}}_{2n }\Big]^{p_{2n}}A^{\frac r 2}_{2n }\Big\}^{w_{2}}$;\\
\indent (II.3) $A^{r-t_{n}}_{2n }\geq     \Big\{A^{\frac r 2}_{2n }\Big[A^{-{\frac {t_{n}}{2}}}_{2n }\big\{A^{\frac {t_{n-1}}{2}}_{2n } \cdots A^{\frac {t_{2}}{2}}_{7}\big[A^{-{\frac {t_{2}}{2}}}_{6}\cdot\{A^{\frac {t_{1}}{2}}_{5}(A^{-{\frac {t_{1}}{2}}}_{4}A^{p_{1}}_{3}A^{-{\frac {t_{1}}{2}}}_{4})^{p_{2}} A^{\frac {t_{1}}{2}}_{5}\}^{p_{3}} \cdot\\
 A^{-{\frac {t_{2}}{2}}}_{6}\big]^{p_{4}}A^{\frac {t_{2}}{2}}_{7}
\cdots A^{\frac {t_{n-1}}{2}}_{2n }\big\}^{p_{2n-1}}A^{-{\frac {t_{n}}{2}}}_{2n }\Big]^{p_{2n}}A^{\frac r 2}_{2n }\Big\}^{w_{3}}$;\\
\indent  {  $\cdots \cdots \cdots \cdots$}\\
\indent (II.n) $A^{r-t_{n}}_{2n }\geq     \Big\{A^{\frac r 2}_{2n }\Big[A^{-{\frac {t_{n}}{2}}}_{2n }\big\{A^{\frac {t_{n-1}}{2}}_{2n } \cdots A^{\frac {t_{2}}{2}}_{n+4}\big[A^{-{\frac {t_{2}}{2}}}_{n+3}\{A^{\frac {t_{1}}{2}}_{n+2}(A^{-{\frac {t_{1}}{2}}}_{n+1}A^{p_{1}}_{n}A^{-{\frac {t_{1}}{2}}}_{n+1})^{p_{2}} A^{\frac {t_{1}}{2}}_{n+2}\}^{p_{3}}  \\
 A^{-{\frac {t_{2}}{2}}}_{n+3}\big]^{p_{4}}A^{\frac {t_{2}}{2}}_{n+4}
\cdots A^{\frac {t_{n-1}}{2}}_{2n }\big\}^{p_{2n-1}}A^{-{\frac {t_{n}}{2}}}_{2n }\Big]^{p_{2n}}A^{\frac r 2}_{2n }\Big\}^{w_{n}}$;\\
\indent (II.n+1) $A^{r-t_{n}}_{ 1}\leq     \Big\{A^{\frac r 2}_{ 1}\Big[A^{-{\frac {t_{n}}{2}}}_{1}\big\{A^{\frac {t_{n-1}}{2}}_{ 1 } \cdots A^{\frac {t_{2}}{2}}_{n-2}\big[A^{-{\frac {t_{2}}{2}}}_{n-1}\{A^{\frac {t_{1}}{2}}_{n }(A^{-{\frac {t_{1}}{2}}}_{n+1}A^{p_{1}}_{n+2}A^{-{\frac {t_{1}}{2}}}_{n+1})^{p_{2}} A^{\frac {t_{1}}{2}}_{n }\}^{p_{3}}  \\
 A^{-{\frac {t_{2}}{2}}}_{n-1}\big]^{p_{4}}A^{\frac {t_{2}}{2}}_{n-2}
\cdots A^{\frac {t_{n-1}}{2}}_{1}\big\}^{p_{2n-1}}A^{-{\frac {t_{n}}{2}}}_{1}\Big]^{p_{2n}}A^{\frac r 2}_{1}\Big\}^{w_{n+1}}$;\\
\indent  {  $\cdots \cdots \cdots \cdots$}\\
\indent (II.2n-2) $A^{r-t_{n}}_{ 1}\leq     \Big\{A^{\frac r 2}_{
1}\Big[A^{-{\frac {t_{n}}{2}}}_{1}\big\{A^{\frac {t_{n-1}}{2}}_{1}
\cdots A^{\frac {t_{2}}{2}}_{2n-5}\big[A^{-{\frac
{t_{2}}{2}}}_{2n-4}\{A^{\frac {t_{1}}{2}}_{2n-3 }\cdot(A^{-{\frac
{t_{1}}{2}}}_{2n-2}A^{p_{1}}_{2n-1}A^{-{\frac
{t_{1}}{2}}}_{2n-2})^{p_{2}}\cdot \\ A^{\frac
{t_{1}}{2}}_{2n-3}\}^{p_{3}}
 A^{-{\frac {t_{2}}{2}}}_{2n-4}\big]^{p_{4}}A^{\frac {t_{2}}{2}}_{2n-5}
\cdots A^{\frac {t_{n-1}}{2}}_{1}\big\}^{p_{2n-1}}A^{-{\frac {t_{n}}{2}}}_{1}\Big]^{p_{2n}}A^{\frac r 2}_{1}\Big\}^{w_{2n-2}}$;\\
\indent (II.2n-1) $A^{r-t_{n}}_{ 1}\leq     \Big\{A^{\frac r 2}_{
1}\Big[A^{-{\frac {t_{n}}{2}}}_{1}\big\{A^{\frac {t_{n-1}}{2}}_{2}
\cdots A^{\frac {t_{2}}{2}}_{2n-4}\big[A^{-{\frac
{t_{2}}{2}}}_{2n-3}\{A^{\frac {t_{1}}{2}}_{2n-2 }\cdot(A^{-{\frac
{t_{1}}{2}}}_{2n-1}A^{p_{1}}_{2n}A^{-{\frac
{t_{1}}{2}}}_{2n-1})^{p_{2}}\cdot \\ A^{\frac
{t_{1}}{2}}_{2n-2}\}^{p_{3}}
 A^{-{\frac {t_{2}}{2}}}_{2n-3}\big]^{p_{4}}A^{\frac {t_{2}}{2}}_{2n-4}
\cdots A^{\frac {t_{n-1}}{2}}_{2}\big\}^{p_{2n-1}}A^{-{\frac {t_{n}}{2}}}_{1}\Big]^{p_{2n}}A^{\frac r 2}_{1}\Big\}^{w_{2n-1}}$,\\
then the operator order $A_{2n}\geq A_{2n-1}\geq A_{2n-2}\cdots \geq
A_{3}\geq A_{2}\geq A_{1}$ holds.

\noindent{\bf Proof.} {Replace $A_{2n+1}$ by $A_{2n}$ } in Theorem 2.1.          $\square$\\
\noindent {\bf Remark 2.2.} If $w_{1}=w_{2}=\cdots=w_{2n-1}={\frac
{r-t_{n}}{\psi [2n]-t_{n}+r}}$, where $\psi
[2n]=\{\cdots[\{[(p_{1}-t_{1})p_{2}+t_{1}]p_{3}-t_{2}\}p_{4}+t_{2}]p_{5}-\cdots
-t_{n}\}p_{2n}+t_{n}$, the condition in Theorem 2.2 is the
sufficient and necessary condition of $A_{2n}\geq A_{2n-1}\geq
A_{2n-2}\geq \cdots\geq A_{3}\geq A_{2}\geq A_{1}$. See \cite{Shi}
for details.

\noindent {\bf Remark 2.3.} Together Theorem 2.1 with Theorem 2.2,
we list the sufficient condition of $A_{k}\geq A_{k-1}\geq
\cdots\geq A_{3}\geq A_{2}\geq A_{1}$ for any integer $k$.

\begin{center}

\end{center}
\end{document}